\documentclass[twoside,10pt]{article}
\usepackage{amssymb}
\usepackage{amsthm}
\usepackage{amsmath}
\usepackage[all]{xy}
\usepackage{enumerate}
\usepackage{fancyhdr}
\usepackage{a4wide}

\newtheorem{theorem}{Theorem}[section]
\newtheorem{proposition}[theorem]{Proposition}
\newtheorem{corollary}[theorem]{Corollary}
\newtheorem{lemma}[theorem]{Lemma}
\theoremstyle{definition}

\newtheorem*{Beweis}{Proof}
\newtheorem{definition}[theorem]{Definition}
\newtheorem{punto}[theorem]{}
\theoremstyle{remark}
\newtheorem{remark}[theorem]{Remark}
\newtheorem{ex}[theorem]{Example}
\newtheorem{exs}[theorem]{Examples}
\newtheorem{remarks}[theorem]{Remarks}
\CompileMatrices
\swapnumbers
\CompileMatrices
\input{tcilatex}
\begin{document}

\title{Exact Sequences of Semimodules over Semirings}
\author{\textbf{Jawad Y. Abuhlail\thanks{%
The author would like to acknowledge the support provided by the Deanship of
Scientific Research (DSR) at King Fahd University of Petroleum $\&$ Minerals
(KFUPM) for funding this work through project No. FT100004.}} \\
Department of Mathematics and Statistics\\
King Fahd University of Petroleum $\&$ Minerals\\
abuhlail@kfupm.edu.sa}
\date{}
\maketitle

\begin{abstract}
In this paper, we introduce and investigate a new notion of exact sequences
of semimodules over semirings relative to the canonical image factorization.
Several homological results are proved using the new notion of exactness
including some restricted versions of the Short Five Lemma and the Snake
Lemma opening the door for introducing and investigating \emph{homology
objects} in such categories. Our results apply in particular to the variety
of commutative monoids extending results in homological varieties.
\end{abstract}

\section*{Introduction}

\qquad Semirings and categories of semimodules over them gained recently
increasing interest due to their role in several emerging areas of research
like Idempotent Analysis (\emph{e.g.} \cite{KM1997}, \cite{LMS2001}, \cite%
{Lit2007}), Tropical Geometry (\emph{e.g.} \cite{R-GST2005}, \cite{Mik2006})
and other aspects of modern Mathematics and Mathematical Physics (\emph{e.g.}
\cite{Go19l99a}, \cite{LM2005}). From the categorical (homological) algebra
point of view, several notions of exact sequences of semimodules were
considered in the literature (\emph{e.g.} \cite{Tak1981}, \cite{Pat2003}, 
\cite{PD2006}). However, none these notions enabled a smooth development of
a homological theory for semimodules over semirings.

In this manuscript, and based on investigations on the notion of exact
sequences in arbitrary non-exact categories w.r.t. a given factorization
system $(\mathbf{E},\mathbf{M})$ \cite{AHS2004}, we provide a new notion of
exactness for semimodules over semirings w.r.t. the canonical image
factorization. We illustrate the usefulness of this new notion by proving
some restricted versions of the Short Five Lemma and the Snake Lemma.

The manuscript is divided as follows. After this brief introduction, and for
the convention of the reader, we collect in Section 1 some definitions and
results on semirings and semimodules and clarify the differences between the
terminology used in this paper and the terminology of \cite{Tak1981} and 
\cite{Go19l99a}; we also clarify the reason for changing some terminology.
In Section 2, we introduce our new notion of exact sequences of semimodules
over semirings. We demonstrate how this notion enables us to characterize in
a very simple way, similar to that in homological categories, different
classes of morphisms (e.g. monomorphisms, regular epimorphisms,
isomorphisms). In Section 3, we illustrate the main advantages of our notion
of exactness over the existing ones by showing how it enables us to prove
some of the elementary diagram lemmas for semimodules over semirings.
Moreover, we introduce a restricted version of the \emph{Short Five Lemma }%
(Proposition \ref{short-5}) which characterizes the homological categories
among the pointed regular ones \cite{BB2004}. Moreover, we prove a
restricted version of the \emph{Snake Lemma }(Theorem \ref{snake}) for
cancellative semimodules (cancellative commutative monoids) which opens the
door for introducing and investigating \emph{homology objects} in such
categories; for (co)homology monoids see for example \cite{Ina1997}, \cite%
{Pat2000a} and \cite{Pat2006}.

\section{Semirings and Semimodules}

\qquad \emph{Semirings} (\emph{semimodules}) are roughly speaking, rings
(modules) without subtraction. Semirings were studied by Vandiver (\emph{e.g.%
} \cite{Van1934}, \cite{Van1936}) as they provide a natural unification
rings and bounded distributive lattices. Since then, semirings were shown to
have significant applications in several areas as Automata Theory (\emph{e.g.%
} \cite{Eil1974}, \cite{Eil1976}, \cite{KS1986}), Theoretical Computer
Science (\emph{e.g.} \cite{HW1998}) and several areas of mathematics (\emph{%
e.g.} \cite{Go19l99a}, \cite{Gol1999b}). Recent applications in emerging
areas of research are demonstrated in several manuscripts (\emph{e.g.} \cite%
{KM1997}, \cite{LMS2001}, \cite{Gol2003}, \cite{LM2005}, \cite{R-GST2005}, 
\cite{Mik2006}, \cite{Lit2007}). Moreover, Durov demonstrated in his
dissertation \cite{Dur2007} that semirings are in one-to-one correspondence
with what he called \emph{algebraic additive monads} on the category $%
\mathbf{Set}$ of sets. A connection between semirings and the so-called $%
\mathbb{F}$-rings, where $\mathbb{F}$ is the field with one element, was
pointed out in \cite[1.3 -- 1.4]{PL2011}.

The basic concepts in the theory of semimodules over semirings were
developed in several articles (\emph{e.g.} \cite{Tak1981}, \cite{Tak1982a}, 
\cite{Tak1982a}, \cite{Tak1982b}, \cite{Tak1983}, \cite{Pat1998}, \cite%
{Pat2000a}, \cite{Pat2000b}, \cite{Pat2003}, \cite{Pat2006}, \cite{Pat2009}, 
\cite{Kat2004b}, \cite{KTN2009}, \cite{KN2011}, \cite{IK2011}, \cite{IK2011}%
). In what follows, we collect and introduce some terminology that will be
needed in the sequel.

\begin{punto}
Let $(S,+)$ be an Abelian additive semigroup. We call $s\in S$ \emph{%
cancellable} iff for any $s_{1},s_{2}\in S$ we have $[s+s_{1}=s+s_{2}%
\Longrightarrow s_{1}=s_{2}].$ We call $S$ \emph{cancellative} iff all
elements of $S$ are cancellable. We say that a morphism of abelian additive
semigroups $f:S\longrightarrow S^{\prime }$ is \emph{cancellative} iff $%
f(s)\in S^{\prime }$ is cancellable for every $s\in S.$ The\emph{\
subtractive closure} of a non-empty subset $X\subseteq S$ is given by%
\begin{equation*}
\overline{X}:=\{s\in S\mid s+x_{1}=x_{2}\text{ for some }x_{1},x_{2}\in X\}.
\end{equation*}%
If $X$ is a subsemigroup of $S,$ then we say that $X$ is \emph{subtractive }%
iff $X=\overline{X}.$ We call a morphism of Abelian semigroups $%
f:S\longrightarrow S^{\prime }$ \emph{subtractive} iff $f(S)\subseteq
S^{\prime }$ is subtractive.
\end{punto}

\begin{punto}
A \emph{semiring} is an algebraic structure $(S,+,\cdot ,0,1)$ consisting of
a non-empty set $S$ with two binary operations \textquotedblleft $+$%
\textquotedblright\ (addition) and \textquotedblleft $\cdot $%
\textquotedblright\ (multiplication) satisfying the following conditions:

\begin{enumerate}
\item $(S,+,0)$ is an Abelian monoid with neutral element $0;$

\item $(S,\cdot ,1)$ is a monoid with neutral element $1;$

\item $x\cdot (y+z)=x\cdot y+x\cdot z$ and $(y+z)\cdot x=y\cdot x+z\cdot x$
for all $x,y,z\in S;$

\item $0\cdot s=0=s\cdot 0$ for every $s\in S$ (\emph{i.e.} $0$ is \emph{%
absorbing}).
\end{enumerate}

Let $S,S^{\prime }$ be semirings. A map $f:S\rightarrow S^{\prime }$ is said
to be a \emph{morphism of semirings} iff for all $s_{1},s_{2}\in S:$%
\begin{equation*}
f(s_{1}+s_{2})=f(s_{1})+f(s_{2}),\text{ }f(s_{1}s_{2})=f(s_{1})f(s_{2}),%
\text{ }f(0_{S})=0_{S^{\prime }}\text{ and }f(1_{S})=1_{S^{\prime }}.
\end{equation*}
\end{punto}

\begin{punto}
Let $(S,+,\cdot )$ be a semiring. We say that $S$ is

\emph{cancellative} iff the additive semigroup $(S,+)$ is cancellative;

\emph{commutative} iff the multiplicative semigroup $(S,\cdot )$ is
commutative;

\emph{semifield} iff $(S\backslash \{0\},\cdot ,1)$ is a commutative group.
\end{punto}

\begin{exs}
Rings are indeed semirings. A trivial example of a \emph{commutative }%
semiring that is not a ring is $(\mathbb{N}_{0},+,\cdot )$ (the set of
non-negative integers). Indeed, $(\mathbb{R}_{0}^{+},+,\cdot )$ and $(%
\mathbb{Q}_{0}^{+},+,\cdot )$ are semifields. A more interesting example is
the semi-ring $(\mathrm{Ideal}(R),+,\cdot )$ consisting of all ideals of a
not necessarily commutative ring (\emph{Dedekind} \cite{Ded1894}). On the
other hand, for an integral domain $R,$ $(\mathrm{Ideal}(R),+,\cap )$ is a
semiring if and only if $R$ is a Pr\"{u}fer domain. Every bounded
distributive lattice $(R,\vee ,\wedge )$ is a commutative (additively
idempotent) semiring. The \emph{additively idempotent} semirings $\mathbb{R}%
_{\max }:=(\mathbb{R}\cup \{-\infty \},\max ,+)$ and $\mathbb{R}_{\min }:=(%
\mathbb{R}\cup \{\infty \},\min ,+)$ play important roles idempotent and
tropical mathematics (\emph{e.g.} \cite{Lit2007}); the subsemirings $\mathbb{%
N}_{\max }:=(\mathbb{N}\cup \{-\infty \},\max ,+)$ and $\mathbb{N}_{\min }:=(%
\mathbb{N}\cup \{\infty \},\min ,+)$ show up in Automata Theory (\emph{e.g.} 
\cite{Eil1974}, \cite{Eil1976}). In the sequel, we always assume that $%
0_{S}\neq 1_{S}$ so that $S\neq \{0\},$ the \emph{zero semiring}.
\end{exs}

\begin{punto}
Let $S$ be a semiring. A \emph{right }$S$\emph{-semimodule} is an algebraic
structure $(M,+,0_{M};\leftharpoondown )$ consisting of a non-empty set $M,$
a binary operation \textquotedblleft $+$\textquotedblright\ along with a
right $S$-action%
\begin{equation*}
M\times S\longrightarrow M,\text{ }(m,s)\mapsto ms,
\end{equation*}%
such that:

\begin{enumerate}
\item $(M,+,0_{M})$ is an Abelian monoid with neutral element $0_{M};$

\item $(ms)s^{\prime }=m(ss^{\prime }),$ $(m+m^{\prime })s=ms+m^{\prime }s$
and $m(s+s^{\prime })=ms+ms^{\prime }$ for all $s,s^{\prime }\in S$ and $%
m,m^{\prime }\in M;$

\item $m1_{S}=m$ and $m0_{S}=0_{M}=0_{M}s$ for all $m\in M$ and $s\in S.$

Let $M,M^{\prime }$ be right $S$-semimodules. A map $f:M\rightarrow
M^{\prime }$ is said to be a \emph{morphism of right }$S$\emph{-semimodules}
(or $S$\emph{-linear}) iff for all $m_{1},m_{2}\in M$ and $s\in S:$%
\begin{equation*}
f(m_{1}+m_{2})=f(m_{1})+f(m_{2})\text{ and }f(ms)=f(m)s.
\end{equation*}%
One can easily check that for any morphism of semimodules $%
f:M\longrightarrow M^{\prime },$ we have $M^{\prime }/\overline{f(M)}%
=M^{\prime }/f(M).$ The set $\mathrm{Hom}_{S}(M,M^{\prime })$ of $S$-linear
maps from $M$ to $M^{\prime }$ is clearly a monoid under addition. The
category of right $S$-semimodules is denoted by $\mathbb{S}_{S}.$ Similarly,
one can define the category of left $S$-semimodules $_{S}\mathbb{S}.$ A
right $S$-semimodule $M_{S}$ is said to be \emph{cancellative} iff the
semigroup $(M,+)$ is cancellative. With $\mathbb{CS}_{S}\subseteq \mathbb{S}%
_{S}$ (resp. $_{S}\mathbb{CS}\subseteq $ $_{S}\mathbb{S})$ we denote the
full subcategory of cancellative right (left) $S$-semimodules.
\end{enumerate}
\end{punto}

\begin{ex}
Every Abelian monoid $(M,+,0_{M})$ is an $\mathbb{N}_{0}$-semimodule in the
obvious way. Moreover, the categories $\mathbf{CMon}$ of commutative monoids
and the category $\mathbb{S}_{\mathbb{N}_{0}}$ of $\mathbb{N}_{0}$%
-semimodules are isomorphic.
\end{ex}

\begin{punto}
Let $M$ be a right $S$-semimodule. A non-empty subset $L\subseteq M$ is said
to be an $S$\emph{-subsemimodule, }and we write $L\leq _{S}M$ iff $L$ is
closed under \textquotedblleft $+_{M}$\textquotedblright\ and $ls\in L$ for
all $l\in L$ and $s\in S.$ Every $S$-subsemimodule $L\leq _{S}M$ induces an $%
S$\emph{-congruence} on $M$ (\emph{e.g.} \cite{Gol1999b}) given by the \emph{%
Bourne relation}%
\begin{equation*}
m_{1}\equiv _{L}m_{2}\Leftrightarrow m_{1}+l_{1}=m_{2}+l_{2}\text{ for some }%
l_{1},l_{2}\in L.
\end{equation*}%
We call the $S$-semimodule $M/L:=M/_{\equiv _{L}}$ the \emph{quotient of }$M$%
\emph{\ by }$L$ or the \emph{factor semimodule} of $M$ by $L.$ If $M$ is
cancellative, then $L$ and $M/L$ are indeed cancellative.
\end{punto}

\begin{proposition}
The category $\mathbb{S}_{S}$ and its full subcategory $\mathbb{CS}_{S}$
have kernels and cokernels, where for any morphism of $S$-semimodules $%
f:M\rightarrow N$ we have%
\begin{equation*}
\mathrm{Ker}(f)=\{m\in M\mid f(m)=0\}\text{ and }\mathrm{\mathrm{Coker}}%
(f)=N/f(M).
\end{equation*}
\end{proposition}

\begin{punto}
The category of $S$-semimodules is regular (\emph{e.g.} \cite{Gri1971}, \cite%
{Bor1994}); in particular, $\mathbb{S}_{S}$ has a $(\mathbf{Surj},\mathbf{%
Mono})$\emph{-factorization structure} \cite{AHS2004}. Let $\gamma
:X\longrightarrow Y$ be a morphism of $S$-semimodules. Then the image and
coimage of $\gamma $ are given respectively by $\func{Im}(\gamma )=\gamma
(X) $ and $\mathrm{Coim}(\gamma )=X/f,$ where $X/f$ is the quotient
semimodule $X/\equiv _{f}$ and $x\equiv _{f}x^{\prime }$ iff $%
f(x)=f(x^{\prime }).$ Indeed, we have a canonical isomorphism%
\begin{equation*}
d_{\gamma }:\mathrm{Coim}(\gamma )\simeq \func{Im}(\gamma ),\text{ }%
[x]\mapsto \gamma (x).
\end{equation*}
\end{punto}

\begin{remark}
\label{image}For any morphisms of $S$-semimodules $\gamma :X\longrightarrow
Y $ we have%
\begin{equation*}
\begin{tabular}{lll}
$\mathrm{Ker}(\mathrm{\mathrm{coker}}(\gamma ))$ & $=$ & $\{y\in Y\mid
y\equiv _{\gamma (X)}0\}$ \\ 
& $=$ & $\{y\in Y\mid y+\gamma (x_{1})=\gamma (x_{2})$ for some $%
x_{1},x_{2}\in X\}$ \\ 
& $=$ & $\overline{\gamma (X)}.$%
\end{tabular}%
\end{equation*}%
Takahashi \cite{Tak1981} defined the \emph{image} of $\gamma $ as $\overline{%
\gamma (X)}.$ While $\mathrm{Ker}(\mathrm{coker}(\gamma ))$ would have been
the \emph{categorical image} (\emph{e.g.} \cite[5.8.7]{Fai1973}, \cite%
{EW1987}) if $\mathbb{S}_{S}$ were a Puppe-exact category (\emph{e.g.} \cite%
{Pup1962}, \cite{Sch1972}), which is indeed not the case in general, the
image of $\gamma $ is in fact the proper image $\gamma (X).$
\end{remark}

\begin{punto}
We call a morphism of $S$-semimodules $\gamma :X\longrightarrow Y:$

$k$\emph{-uniform} iff for any $x_{1},x_{2}\in X:$%
\begin{equation}
\gamma (x_{1})=\gamma (x_{2})\Longrightarrow \text{ }\exists \text{ }%
k_{1},k_{2}\in \mathrm{Ker}(\gamma )\text{ s.t. }x_{1}+k_{1}=x_{2}+k_{2};
\label{k-steady}
\end{equation}

$i$\emph{-uniform} iff $\gamma (X)=\overline{\gamma (X)};$

\emph{uniform }iff $\gamma $ is $k$-uniform and $i$-uniform;

\emph{semi-monomorphism} iff $\mathrm{Ker}(\gamma )=0;$

\emph{semi-epimorphism} iff $\overline{\gamma (X)}=Y;$

\emph{semi-isomorphism} iff $\mathrm{Ker}(\gamma )=0$ and$\ \overline{\gamma
(X)}=Y.$
\end{punto}

\qquad The following lemma provides a description of the above mentioned
classes of morphisms in terms of the well-known classes of (normal, regular)
monomorphisms and epimorphisms. The proof follows immediately from Remark %
\ref{image} and the fact that the canonical $(\mathbf{Surj},\mathbf{Mono})$-
factorization of $\gamma $ is given by $\gamma :X\overset{\mathrm{coim}%
(\gamma )}{\longrightarrow }\gamma (X)\overset{\mathrm{im}(\gamma )}{%
\hookrightarrow }Y.$

\begin{lemma}
\label{co-stead}Let $\gamma :X\longrightarrow Y$ be a morphism of $S$%
-semimodules.

\begin{enumerate}
\item The following are equivalent:

\begin{enumerate}
\item $\gamma =m\circ \mathrm{coker}(\mathrm{ker}(\gamma ))$ is the $(%
\mathbf{Surj},\mathbf{Mono})$- factorization of $\gamma ,$ where $m$ is an
appropriate monomorphism;

\item $\mathrm{Coker}(\mathrm{ker}(\gamma ))\simeq \mathrm{Coim}(\gamma );$

\item $X/\mathrm{Ker}(\gamma )\simeq \gamma (X);$

\item $\gamma $ is $k$-uniform.
\end{enumerate}

\item The following are equivalent:

\begin{enumerate}
\item $\gamma =\mathrm{ker}(\mathrm{coker}(\gamma ))\circ e$ is the $(%
\mathbf{Surj},\mathbf{Mono})$- factorization of $\gamma ,$ where $e$ is an
appropriate regular epimorphism;

\item $\mathrm{Ker}(\mathrm{coker}(\gamma ))\simeq \func{Im}(\gamma );$

\item $\overline{\gamma (X)}=\gamma (X);$

\item $\gamma $ is $i$-uniform (subtractive).
\end{enumerate}

\item The following are equivalent:

\begin{enumerate}
\item $\gamma =m\circ \mathrm{coker}(\mathrm{ker}(\gamma ))=\mathrm{ker}(%
\mathrm{coker}(\gamma ))\circ e$ is the $(\mathbf{Surj},\mathbf{Mono})$-
factorization of $\gamma $ for an appropriate monomorphism $m$ and an
appropriate regular epimorphism $e;$

\item $\mathrm{Coker}(\mathrm{ker}(\gamma ))\simeq \mathrm{Ker}(\mathrm{Coker%
}(\gamma ));$

\item $X/\mathrm{Ker}(\gamma )\simeq \overline{\gamma (X)};$

\item $\gamma $ is uniform.
\end{enumerate}
\end{enumerate}
\end{lemma}

\begin{remarks}
\begin{enumerate}
\item Lemma \ref{co-stead} describes a $k$-uniform ($i$-uniform) morphism of
semimodules as a composition of a normal epimorphism followed by a
monomorphism (a regular epimorphism followed by a normal monomorphism). We
use this terminology because it is brief and related to the usual
terminology used in the literature of semimodules (see \textquotedblleft
2\textquotedblright\ below).

\item The uniform ($k$-uniform, $i$-uniform) morphisms of semimodules were
called \emph{regular} ($k$\emph{-regular, }$i$-\emph{regular}) by Takahashi 
\cite{Tak1982c}. We think that our terminology avoids confusion since a
regular monomorphism (regular epimorphism) has a different well-established
meaning in the language of Category Theory (\emph{e.g.} \cite{AHS2004}).
\end{enumerate}
\end{remarks}

\begin{punto}
Let $M$ be an $S$-semimodule, $L\leq _{S}M$ an $S$-subsemimodule and
consider the factor semimodule $M/L.$ Then we have a surjective morphism of $%
S$-semimodules%
\begin{equation*}
\pi _{L}:=M\rightarrow M/L,\text{ }m\mapsto \lbrack m]
\end{equation*}%
with%
\begin{equation*}
\mathrm{Ker}(\pi _{L})=\{m\in M\mid m+l_{1}=l_{2}\text{ for some }%
l_{1},l_{2}\in L\}=\overline{L};
\end{equation*}%
in particular, $L=\mathrm{Ker}(\pi _{L})$ if and only if $L\subseteq M$ is
subtractive.
\end{punto}

\section{Exact Sequences of Semimodules}

\qquad Throughout this section, $S$ is a semiring, and an $S$-semimodule is
a right $S$-semimodule unless otherwise explicitly specified. Moreover, $%
\mathbb{S}_{S}$ ($\mathbb{CS}_{S}$) denotes the category of (cancellative)
right $S$-semimodules. For undefined terminology from Category Theory, we
refer to \cite{Mac1998} and \cite{AHS2004}.

The notion of \emph{exact sequences} of semimodules adopted by Takahashi 
\cite{Tak1981} ($L\overset{f}{\longrightarrow }M\overset{g}{\longrightarrow }%
N$ is exact iff $\overline{f(M)}=\mathrm{Ker}(g)$) seems to be inspired by
the definition of exact sequences in Puppe-exact categories \cite{Pup1962}
(see also \cite{Bue2010}). We believe it is inappropriate. The reason for
this is that neither $\mathrm{Ker}(\mathrm{coker}(f))=\overline{f(L)}$ is
the appropriate \emph{image} of $f$ nor is $\mathrm{Coker}(\mathrm{ker}%
(g))=B/\mathrm{Ker}(g)$ the appropriate \emph{coimage} of $g.$

Being a Barr-exact category, a natural tool to study exactness in the
category of semimodules is that of an \emph{exact fork} \cite{Bar1971} and
applied to study exact functors between categories of semimodules (\emph{e.g.%
} Katsov et al. \cite{KN2011}). However, since the category of semimodules
has additional features, one still expects to deal with exact sequences
rather than the more complicated exact forks.

In addition to Takahashi's definition of exact sequences of semimodules, two
different notions of exactness for sequences of semimodules over semirings
are available in the literature: one notion is used by Patchkoria \cite%
{Pat2003} ($L\overset{f}{\longrightarrow }M\overset{g}{\longrightarrow }N$
is exact iff $f(L)=\mathrm{Ker}(g)$) and another can be found in \cite%
{PD2006} ($L\overset{f}{\longrightarrow }M\overset{g}{\longrightarrow }N$ is
exact iff ${\overline{f(L)}}=\mathrm{Ker}(g)$ and $g$ is $k$-uniform). Each
of these definitions is stronger than Takahashi's notion of exactness and
each proved to be more efficient in establishing some homological results
for semimodules over semirings. However, no clear \emph{categorical}
justification for choosing either of these two definitions was provided. A
closer look at these definitions shows that they are in fact dual to each
other in some sense, and so it not logical to choose one of them and drop
the other. This motivated us to introduce a new notion of exact sequences of
semimodules as a combination of the two above mentioned notions of exact
sequences of semimodules. Our notion is motivated by intensive
investigations of exact sequences in arbitrary not necessarily Puppe-exact
pointed categories. For the record, we mention here that there are some
other notions of exact sequences of semimodules in the literature which we
did not mention here since their definitions are very technical and not
categorical.

\begin{definition}
We say that a sequence $X\overset{f}{\longrightarrow }Y\overset{g}{%
\longrightarrow }Z$ of semimodules is \emph{exact} (more precisely $(\mathbf{%
Surj},\mathbf{Inj})$\emph{-exact}) iff $f(X)=\mathrm{Ker}(g)$ and $g$ is $k$%
-uniform.
\end{definition}

\begin{lemma}
\label{S-mon-epi}Let%
\begin{equation}
L\overset{f}{\rightarrow }M\overset{g}{\rightarrow }N  \label{lmn}
\end{equation}%
be a sequence of $S$-semimodules with $g\circ f=0$ and consider the induced
morphisms $f^{\prime }:L\rightarrow \mathrm{Ker}(g)$ and $g^{\prime \prime }:%
\mathrm{Coker}(f)\rightarrow N.$

\begin{enumerate}
\item If $f^{\prime }$ is an epimorphism, then $\overline{f(L)}=\mathrm{Ker}%
(g).$

\item $f^{\prime }$ is a regular epimorphism (surjective) if and only if $%
f(L)=\mathrm{Ker}(g)$ if and only if $\overline{f(L)}=\mathrm{Ker}(g)$ and $%
f $ is $i$-uniform.

\item $g^{\prime \prime }:\mathrm{Coker}(f)\rightarrow N$ is a monomorphism
if and only if $\overline{f(L)}=\mathrm{Ker}(g)$ and $g$ is $k$-uniform.
\end{enumerate}
\end{lemma}

\begin{Beweis}
Since $g\circ f=0,$ we have $f(L)\subseteq \overline{f(L)}\subseteq \mathrm{%
Ker}(g).$

\begin{enumerate}
\item Assume that $f^{\prime }:L\rightarrow \mathrm{Ker}(g)$ is an
epimorphism. Suppose that $\overline{f(L)}\subsetneqq \mathrm{Ker}(g),$ so
that there exists $m^{\prime }\in \mathrm{Ker}(g)\backslash \overline{f(L)}.$
Consider the $S$-linear maps 
\begin{equation*}
L\overset{f^{\prime }}{\rightarrow }\mathrm{Ker}(g)\overset{f_{1}}{\underset{%
f_{2}}{\rightrightarrows }}\mathrm{Ker}(g)/f(L),
\end{equation*}%
where $f_{1}(m)=[m]$ and $f_{2}(m)=[0]$ for all $m\in \mathrm{Ker}(g).$ For
each $l\in L$ we have 
\begin{equation*}
(f_{1}\circ f^{\prime })(l)=[f(l)]=[0]=(f_{2}\circ f^{\prime })(l).
\end{equation*}%
Whence, $f_{1}\circ f^{\prime }=f_{2}\circ f^{\prime }$ while $f_{1}\neq
f_{2}$ (since $f_{1}(m^{\prime })=[m^{\prime }]\neq \lbrack
0]=f_{2}(m^{\prime });$ otherwise $m^{\prime }+f(l_{1})=f(l_{2})$ for some $%
l_{1},l_{1}^{\prime }\in L$ and $m^{\prime }\in \overline{f(L)}$ which
contradicts our assumption). So, $f^{\prime }$ is not an epimorphism, a
contradiction. Consequently, $\overline{f(L)}=\mathrm{Ker}(g).$

\item Clear.

\item $(\Rightarrow )$ Assume that $g^{\prime \prime }:\mathrm{Coker}%
(f)\rightarrow N$ is a monomorphism. Let $m\in \mathrm{Ker}(g),$ so that $%
g(m)=0.$ Then $g^{\prime \prime }([m])=0.$ Since $g^{\prime \prime }$ is a
monomorphism, we have $[m]=[0]$ and so $m+f(l)=f(l^{\prime })$ for some $%
l,l^{\prime }\in L,$ \emph{i.e.} $m\in \overline{f(L)}.$ Suppose now that $%
g(m)=g(m^{\prime })$ for some $m,m^{\prime }\in M.$ Then $g^{\prime \prime
}([m])=g^{\prime \prime }([m^{\prime }])$ and it follows, by the injectivity
of $g^{\prime \prime },$ that $[m]=[m^{\prime }]$ which implies that $%
m_{1}+m_{1}=m^{\prime }+m_{1}^{\prime }$ for some $m_{1},m_{1}^{\prime }\in 
\overline{f(L)}\subseteq \mathrm{Ker}(g).$ So, $g$ is $k$-uniform.

$(\Leftarrow )$ Assume that $\overline{f(L)}=\mathrm{Ker}(g)$ and that $g$
is $k$-uniform. Suppose that $g^{\prime \prime }([m])=g^{\prime \prime
}([m^{\prime \prime }]),$ whence $g(m)=g(m^{\prime }),$ for some $%
m_{1},m_{2}\in M.$ Since $g$ is $k$-uniform, we have $m+k=m^{\prime
}+k^{\prime }$ for some $k,k^{\prime }\in \mathrm{Ker}(g)=\overline{f(L)}$
and it follows that $[m]=[m^{\prime }].\blacksquare $
\end{enumerate}
\end{Beweis}

\begin{remarks}
\label{CS-mon-epi}

\begin{enumerate}
\item A morphism of cancellative semimodules $h:X\rightarrow Y$ is an
epimorphism in $\mathbb{CS}_{S}$ if and only if $\overline{h(X)}=Y.$ Indeed,
if $h$ is an epimorphism, then it follows by Lemma \ref{S-mon-epi} that $%
\overline{h(X)}=Y$ (take $g:Y\rightarrow 0$ as the zero-morphism). On the
other hand, assume that $\overline{h(L)}=Y.$ Let $Z$ be any cancellative
semimodule and consider any $S$-linear maps 
\begin{equation*}
X\overset{h}{\rightarrow }Y\underset{h_{2}}{\overset{h_{1}}{%
\rightrightarrows }}Z
\end{equation*}%
with $h_{1}\circ h=h_{2}\circ h.$ Let $y\in Y$ be arbitrary. By assumption, $%
y+h(x_{1})=h(x_{2})$ for some $x_{1},x_{2}\in X,$ whence%
\begin{equation*}
h_{1}(y)+(h_{1}\circ h)(x_{1})=(h_{1}\circ h)(x_{2})=(h_{2}\circ
h)(x_{2})=h_{2}(y)+(h_{2}\circ h)(x_{1}).
\end{equation*}%
Since $Z$ is cancellative, we conclude that $h_{1}(y)=h_{2}(y).$

\item Consider the embedding $\iota :\mathbb{N}_{0}\hookrightarrow \mathbb{Z}
$ in $\mathbb{CS}_{\mathbb{N}_{0}}.$ Indeed, $\overline{\mathbb{N}_{0}}=%
\mathbb{Z},$ whence $\iota $ is an epimorphism which is not regular. This
shows that not all epimorphisms of semimodules are surjective \cite{TW1989}.

\item Let $L\overset{f}{\rightarrow }M\overset{g}{\rightarrow }N$ be a
sequence in $\mathbb{CS}_{S}$ with $g\circ f=0.$ By \textquotedblleft
1\textquotedblright , the induced morphism $f^{\prime }:L\rightarrow \mathrm{%
Ker}(g)$ is an epimorphism if and only if $\overline{f(L)}=\mathrm{Ker}(g).$
\end{enumerate}
\end{remarks}

\begin{punto}
\label{def-exact}We call a sequence of $S$-semimodules $L\overset{f}{%
\rightarrow }M\overset{g}{\rightarrow }N:$

\emph{proper-exact} iff $f(L)=\mathrm{Ker}(g);$

\emph{semi-exact} iff $\overline{f(L)}=\mathrm{Ker}(g);$

\emph{uniform} (resp. $k$\emph{-uniform}, $i$\emph{-uniform}) iff $f$ and $g$
are uniform (resp. $k$-uniform, $i$-uniform).
\end{punto}

\begin{punto}
We call a (possibly infinite) sequence of $S$-semimodules 
\begin{equation}
\cdots \rightarrow M_{i-1}\overset{f_{i-1}}{\rightarrow }M_{i}\overset{f_{i}}%
{\rightarrow }M_{i+1}\overset{f_{i+1}}{\rightarrow }M_{i+2}\rightarrow \cdots
\label{chain}
\end{equation}

\emph{chain complex} iff $f_{j+1}\circ f_{j}=0$ for every $j;$

\emph{exact} (resp. \emph{proper-exact}, \emph{semi-exact}) iff each partial
sequence with three terms $M_{j}\overset{f_{j}}{\rightarrow }M_{j+1}\overset{%
f_{j+1}}{\rightarrow }M_{j+2}$ is exact (resp. proper-exact, semi-exact);

\emph{uniform }(resp. $k$-\emph{uniform}, $i$-\emph{uniform}) iff $f_{j}$ is
uniform (resp. $k$-uniform, $i$-uniform) for every $j.$
\end{punto}

\begin{definition}
Let $M$ be an $S$-semimodule.

\begin{enumerate}
\item A subsemimodule $L\leq _{S}M$ is said to be a \emph{uniform }(\emph{%
normal})\emph{\ }$S$-\emph{subsemimodule} iff the embedding $%
0\longrightarrow L\overset{\iota }{\rightarrow }M$ is uniform (normal).

\item A quotient $M/\rho ,$ where $\rho $ is an $S$-congruence relation on $%
M,$ is said to be a \emph{uniform }(\emph{conormal})\emph{\ quotient} iff
the surjection $\pi _{L}:M\rightarrow M/\rho $ is uniform (conormal).
\end{enumerate}
\end{definition}

\begin{remark}
Every normal subsemimodule (conormal quotient) is uniform.
\end{remark}

\qquad The following result can be easily verified.

\begin{lemma}
\label{i-uniform}Let $L\overset{f}{\rightarrow }M\overset{g}{\rightarrow }N$
be a sequence of semimodules.

\begin{enumerate}
\item Let $g$ be injective.

\begin{enumerate}
\item $f$ is $k$-uniform if and only if $g\circ f$ is $k$-uniform.

\item If $g\circ f$ is $i$-uniform (uniform), then $f$ is $i$-uniform
(uniform).

\item Assume that $g$ is $i$-uniform. Then $f$ is $i$-uniform (uniform) if
and only if $g\circ f$ is $i$-uniform (uniform).
\end{enumerate}

\item Let $f$ be surjective.

\begin{enumerate}
\item $g$ is $i$-uniform if and only if $g\circ f$ is $i$-uniform.

\item If $g\circ f$ is $k$-uniform (uniform), then $g$ is $k$-uniform
(uniform).

\item Assume that $f$ is $k$-uniform. Then $g$ is $k$-uniform (uniform) if
and only if $g\circ f$ is $k$-uniform (uniform).
\end{enumerate}
\end{enumerate}
\end{lemma}

\begin{remark}
Let $L\leq _{S}M\leq _{S}N$ be $S$-semimodules. It follows directly from the
previous lemma that if $L$ is uniform in $N,$ then $L$ is a uniform in $M$
as well. Moreover, if $M$ is uniform in $N,$ then $L$ is uniform in $N$ if
and only if $L$ is uniform in $M.$
\end{remark}

Our notion of exactness allows characterization of special classes of
morphisms in a way similar to that in homological categories (compare with 
\cite[Proposition 4.1.9]{BB2004}, \cite[Propositions 4.4, 4.6]{Tak1981}, 
\cite[Proposition 15.15]{Go19l99a}):

\begin{proposition}
\label{inj-surj}Consider a sequence of semimodules%
\begin{equation*}
0\longrightarrow L\overset{f}{\longrightarrow }M\overset{g}{\longrightarrow }%
N\longrightarrow 0.
\end{equation*}

\begin{enumerate}
\item The following are equivalent:

\begin{enumerate}
\item $0\longrightarrow L\overset{f}{\rightarrow }M$ is exact;

\item $\mathrm{Ker}(f)=0$ and $f$ is $k$-uniform;

\item $f$ is injective;

\item $f$ is a monomorphism.
\end{enumerate}

\item $0\longrightarrow L\overset{f}{\longrightarrow }M\overset{g}{%
\longrightarrow }N$ is semi-exact and $f$ is uniform if and only if $L\simeq 
\mathrm{Ker}(g).$

\item $0\longrightarrow L\overset{f}{\longrightarrow }M\overset{g}{%
\longrightarrow }N$ is exact if and only if $L\simeq \mathrm{Ker}(g)$ and $g$
is $k$-uniform.

\item The following are equivalent:

\begin{enumerate}
\item $M\overset{\gamma }{\rightarrow }N\rightarrow 0$ is exact;

\item $\mathrm{Coker}(\gamma )=0$\ and $\gamma $ is $i$-uniform;

\item $\gamma $ is surjective;

\item $\gamma $ is a regular epimorphism;

\item $\gamma $ is a subtractive epimorphism
\end{enumerate}

\item $L\overset{f}{\rightarrow }M\overset{g}{\rightarrow }N\rightarrow 0$
is semi-exact and $g$ is uniform if and only if $N\simeq \mathrm{Coker}(f).$

\item $L\overset{f}{\longrightarrow }M\overset{g}{\longrightarrow }%
N\longrightarrow 0$ is exact if and only if $N\simeq \mathrm{Coker}(f)$ and $%
f$ is $i$-uniform.
\end{enumerate}
\end{proposition}

\begin{corollary}
\label{consist}The following are equivalent:

\begin{enumerate}
\item $0\rightarrow L\overset{f}{\rightarrow }M\overset{g}{\rightarrow }%
N\rightarrow 0$ is a exact sequence of $S$-semimodules;

\item $L\simeq \mathrm{Ker}(g)$ and $\mathrm{\mathrm{\mathrm{\mathrm{Coker}}}%
}(f)\simeq N;$

\item $f$ is injective, $f(L)=\mathrm{Ker}(g),$ $g$ is surjective and ($k$%
-)uniform.

In this case, $f$ and $g$ are uniform morphisms.
\end{enumerate}
\end{corollary}

\begin{remarks}
\begin{enumerate}
\item The equivalence \textquotedblleft 1\textquotedblright\ $%
\Leftrightarrow $ \textquotedblleft 2\textquotedblright\ in Corollary \ref%
{consist} shows that our notion of \emph{short exact sequences} of
semimodules is consistent with that in arbitrary pointed categories in the
sense of \cite[Definition 4.1.5]{BB2004}.

\item Takahashi called a short exact sequence $0\longrightarrow L\overset{f}{%
\longrightarrow }M\overset{g}{\longrightarrow }N\longrightarrow 0$ in our
sense \emph{regular exact }or an \emph{extension} of $N$ by $L$ \cite%
{Tak1982b} (see also \cite{Tak1983}, \cite{Pat2003}).

\item A morphism of semimodules $\gamma :X\longrightarrow Y$ is an
isomorphism if and only if $0\longrightarrow X\longrightarrow
Y\longrightarrow 0$ is exact if and only if $\gamma $ is a uniform
bimorphism. The assumption on $\gamma $ to be uniform cannot be removed
here. For example, the embedding $\iota :\mathbb{N}_{0}\longrightarrow 
\mathbb{Z}$ is a bimorphism of commutative monoids ($\mathbb{N}_{0}$%
-semimodules) which is not an isomorphism. Notice that $\iota $ is \emph{not}
$i$-uniform; in fact $\overline{\iota (\mathbb{N}_{0}})=\mathbb{Z}.$
\end{enumerate}
\end{remarks}

\begin{lemma}
\label{1st-IT}\emph{(Compare with \cite[Proposition 4.3.]{Tak1981})} Let $%
\gamma :X\rightarrow Y$ be a morphism of $S$-semimodules.

\begin{enumerate}
\item The sequence%
\begin{equation}
0\rightarrow \mathrm{Ker}(\gamma )\overset{\mathrm{\ker }(\gamma )}{%
\longrightarrow }X\overset{\gamma }{\rightarrow }Y\overset{\mathrm{\mathrm{%
coker}}(\gamma )}{\longrightarrow }\mathrm{\mathrm{\mathrm{\mathrm{\mathrm{%
Coker}}}}}(\gamma )\rightarrow 0  \label{ker-coker}
\end{equation}

is semi-exact. Moreover, (\ref{ker-coker}) is exact if and only if $\gamma $
is uniform.

\item We have two exact sequences%
\begin{equation*}
0\rightarrow \overline{\gamma (X)}\overset{\mathrm{ker}(\mathrm{\mathrm{coker%
}}(\gamma ))}{\longrightarrow }Y\overset{\mathrm{\mathrm{coker}}(\gamma )}{%
\longrightarrow }Y/\gamma (X)\rightarrow 0.
\end{equation*}%
and%
\begin{equation*}
0\rightarrow \mathrm{Ker}(\gamma )\overset{\mathrm{ker}(\gamma )}{%
\longrightarrow }X\overset{\mathrm{\mathrm{coker}}(\mathrm{ker}(\gamma ))}{%
\longrightarrow }X/\mathrm{Ker}(\gamma )\rightarrow 0.
\end{equation*}
\end{enumerate}
\end{lemma}

\begin{corollary}
\label{reg-sub}\emph{(Compare with \cite[Proposition 4.8.]{Tak1981}) }Let $M$
be an $S$-semimodule.

\begin{enumerate}
\item Let $\rho $ be an $S$-congruence relation on $M$ and consider the
sequence of $S$-semimodules%
\begin{equation*}
0\longrightarrow \mathrm{Ker}(\pi _{\rho })\overset{\iota _{\rho }}{%
\longrightarrow }M\overset{\rho }{\longrightarrow }M/\rho \longrightarrow 0.
\end{equation*}

\begin{enumerate}
\item $0\rightarrow \mathrm{Ker}(\pi _{\rho })\overset{\iota _{\rho }}{%
\longrightarrow }M\overset{\pi _{\rho }}{\longrightarrow }M/\rho \rightarrow
0$ is exact.

\item $M/\rho =\mathrm{Coker}(\iota _{\rho }),$ whence $M/\rho $ is a normal
quotient.
\end{enumerate}

\item Let $L\leq _{S}M$ be an $S$-subsemimodule.

\begin{enumerate}
\item The sequence $0\rightarrow L\overset{\iota }{\longrightarrow }M\overset%
{\pi _{L}}{\longrightarrow }M/L\rightarrow 0$ is semi-exact.

\item $0\rightarrow \overline{L}\overset{\iota }{\longrightarrow }M\overset{%
\pi _{L}}{\longrightarrow }M/L\rightarrow 0$ is exact.

\item The following are equivalent:

\begin{enumerate}
\item $0\rightarrow L\overset{\iota }{\longrightarrow }M\overset{\pi _{L}}{%
\longrightarrow }M/L\rightarrow 0$ is exact;

\item $L\simeq \mathrm{Ker}(\pi _{L});$

\item $0\longrightarrow L\overset{\iota }{\longrightarrow }\overline{L}%
\longrightarrow 0$ is exact;

\item $L$ is a uniform subsemimodule;

\item $L$ is a normal subsemimodule.
\end{enumerate}
\end{enumerate}
\end{enumerate}
\end{corollary}

\section{Homological lemmas}

\qquad In this section we prove some diagram lemmas for semimodules over
semirings. These apply in particular to commutative monoids considered as
semimodules over the semiring of non-negative integers.

Recall that a sequence $A\overset{f}{\longrightarrow }B\overset{g}{%
\longrightarrow }C$ of semimodules is exact iff $f(A)=\mathrm{Ker}(g)$ and $%
g $ is $k$-uniform (equivalently, $f(A)=\mathrm{Ker}(g)$ and $%
g(b)=g(b^{\prime })\Longrightarrow b+f(a)=b^{\prime }+f(a^{\prime })$ for
some $a,a^{\prime }\in A$). In diagram chasing, we sometimes do not mention
where some elements belong if this is clear from the context.

\subsection*{The Five Lemma}

\qquad The following result can be easily proved using \emph{diagram chasing}
(compare \textquotedblleft 2\textquotedblright\ with \cite[Lemma 1.9]%
{Pat2006}).

\begin{lemma}
\label{short}Consider the following commutative diagram of semimodules%
\begin{equation*}
\xymatrix{ & & 0 \ar[d] \\ L_1 \ar[r]^{f_1} \ar[d]_{\alpha_1} & M_1
\ar[r]^{g_1} \ar[d]_{\alpha_2} & N_1 \ar[d]_{\alpha_3} \\ L_2 \ar[r]^{f_2}
\ar[d] & M_2 \ar[r]^{g_2} & N_2 \\ 0 & & }
\end{equation*}%
and assume that the first and the third columns are exact (\emph{i.e.} $%
\alpha _{1}$ is surjective and $\alpha _{3}$ is injective).

\begin{enumerate}
\item Let $\alpha _{2}$ be surjective. If the first row is exact, then the
second row is exact.

\item Let $\alpha _{2}$ be injective. If the second row is exact, then the
first row is exact.

\item Let $a_{2}$ be an isomorphism. The first row is exact if and only if
the second row is exact.
\end{enumerate}
\end{lemma}

\begin{lemma}
\label{diagram}Consider the following commutative diagram of semimodules
with exact rows%
\begin{equation*}
\xymatrix{L_1 \ar[r]^{f_1} \ar[d]_{\alpha_1} & M_1 \ar[r]^{g_1}
\ar[d]_{\alpha_2} & N_1 \ar[d]_{\alpha_3} \\ L_2 \ar[r]^{f_2} & M_2
\ar[r]^{g_2} & N_2}
\end{equation*}

\begin{enumerate}
\item We have:

\begin{enumerate}
\item Let $g_{1}$ and $\alpha _{1}$ be surjective. If $\alpha _{2}$ is
injective, then $\alpha _{3}$ is injective.

\item Let $f_{2}$ be injective and $\alpha _{3}$ be a semi-monomorphism. If $%
\alpha _{2}$ is surjective, then $\alpha _{1}$ is surjective.
\end{enumerate}

\item Let $f_{2}$ be a semi-monomorphism.

\begin{enumerate}
\item If $\alpha _{1}$ and $\alpha _{3}$ are semi-monomorphisms, then $%
\alpha _{2}$ is a semi-monomorphism.

\item Let $f_{1},$ $\alpha _{2}$ be cancellative. If $\alpha _{1},$ $\alpha
_{3}$ and $f_{2}$ are injective, then $\alpha _{2}$ is injective.
\end{enumerate}

\item If $\alpha _{1},$ $\alpha _{3},$ $g_{1}$ are surjective (and $\alpha
_{2}$ is $i$-uniform), then $\alpha _{2}$ is a semi-epimorphism (surjective).
\end{enumerate}
\end{lemma}

\begin{Beweis}
\begin{enumerate}
\item Consider the given commutative diagram.

\begin{enumerate}
\item We claim that $\alpha _{3}$ is injective.

Suppose that $\alpha _{3}(n_{1})=\alpha _{3}(n_{1}^{\prime })$ for some $%
n_{1},n_{1}^{\prime }\in N_{1}.$ Since $g_{1}$ is surjective, $%
n_{1}=g_{1}(m_{1})$ and $n_{1}^{\prime }=g_{1}(m_{1}^{\prime })$ for some $%
m_{1},m_{1}^{\prime }\in M_{1}.$ It follows that $(g_{2}\circ \alpha
_{2})(m_{1})=(g_{2}\circ \alpha _{2})(m_{1}^{\prime }).$ Since $g_{2}$ is $k$%
-uniform and $f_{2}(L_{2})=\mathrm{Ker}(g_{2}),$ there exist $%
l_{2},l_{2}^{\prime }\in L_{2}$ such that $\alpha
_{2}(m_{1})+f_{2}(l_{2})=\alpha _{2}(m_{1}^{\prime })+f_{2}(l_{2}^{\prime
}). $ By assumption, $\alpha _{1}$ is surjective and so there exist $%
l_{1},l_{1}^{\prime }\in L_{1}$ such that $\alpha _{1}(l_{1})=l_{2}$ and $%
\alpha _{1}(l_{1}^{\prime })=l_{2}^{\prime }.$ It follows that%
\begin{equation*}
\begin{array}{rclc}
\alpha _{2}(m_{1})+(f_{2}\circ \alpha _{1})(l_{1}) & = & \alpha
_{2}(m_{1}^{\prime })+(f_{2}\circ \alpha _{1})(l_{1}^{\prime }) &  \\ 
\alpha _{2}(m_{1})+(\alpha _{2}\circ f_{1})(l_{1}) & = & \alpha
_{2}(m_{1}^{\prime })+(\alpha _{2}\circ f_{1})(l_{1}^{\prime }) &  \\ 
m_{1}+f_{1}(l_{1}) & = & m_{1}^{\prime }+f_{1}(l_{1}^{\prime }) & \text{(}%
\alpha _{2}\text{ is injective)} \\ 
g_{1}(m_{1}) & = & g_{1}(m_{1}) & \text{(}g_{1}\circ f_{1}=0\text{)} \\ 
n_{1} & = & n_{1}^{\prime } & 
\end{array}%
\end{equation*}

\item We claim that $\alpha _{1}$ is surjective.

Let $l_{2}\in L_{2}.$ Since $\alpha _{2}$ is surjective, there exists $%
m_{1}\in M_{1}$ such that $f_{2}(l_{2})=\alpha _{2}(m_{1}).$ It follows that 
$0=(g_{2}\circ f_{2})(l_{2})=(g_{2}\circ \alpha _{2})(m_{1})=(\alpha
_{3}\circ g_{1})(m_{1}),$ whence $g_{1}(m_{1})=0$ (since $\alpha _{3}$ is a
semi-monomorphism). Since the first row is exact, $m_{1}=f_{1}(l_{1})$ for
some $l_{1}\in L_{1}$ and so $f_{2}(l_{2})=\alpha _{2}(m_{1})=(\alpha
_{2}\circ f_{1})(l_{1})=(f_{2}\circ \alpha _{1})(l_{1}).$ Since $f_{2}$ is
injective, we have $l_{2}=\alpha _{1}(l_{1}).$
\end{enumerate}

\item Let $f_{2}$ be a semi-monomorphism.

\begin{enumerate}
\item We claim that $\alpha _{2}$ is a semi-monomorphism.

Suppose that $\alpha _{2}(m_{1})=0$ for some $m_{1}\in M_{1}.$ Then $(\alpha
_{3}\circ g_{1})(m_{1})=(g_{2}\circ \alpha _{2})(m_{1})=0,$ whence $%
g_{1}(m_{1})=0$ since $\mathrm{Ker}(\alpha _{3})=0.$ Since the first row is
exact, $m_{1}=f_{1}(l_{1})$ for some $l_{1}\in L_{1}.$ So, $0=\alpha
_{2}(m_{1})=(\alpha _{2}\circ f_{1})(l_{1})=(f_{2}\circ \alpha _{1})(l_{1}),$
whence $l_{1}=0$ since both $f_{2}$ and $\alpha _{1}$ are
semi-monomorphisms; consequently, $m_{1}=f_{1}(l_{1})=0.$

\item We claim that $\alpha _{2}$ is injective.

Suppose that $\alpha _{2}(m_{1})=\alpha _{2}(m_{1}^{\prime })$ for some $%
m_{1},m_{1}^{\prime }\in M_{1}.$ Then $(\alpha _{3}\circ
g_{1})(m_{1})=(g_{2}\circ \alpha _{2})(m_{1})=(g_{2}\circ \alpha
_{2})(m_{1}^{\prime })=(\alpha _{3}\circ g_{1})(m_{1}^{\prime }),$ whence $%
g_{1}(m_{1})=g_{1}(m_{1}^{\prime })$ since $\alpha _{3}$ is injective. Since 
$g_{1}$ is $k$-uniform and $\mathrm{Ker}(g_{1})=f_{1}(L_{1}),$ there exist $%
l_{1},l_{1}^{\prime }\in L_{1}$ such that $m_{1}+f_{1}(l_{1})=m_{1}^{\prime
}+f_{1}(l_{1}^{\prime }).$ Then we have%
\begin{equation*}
\begin{array}{rclc}
\alpha _{2}(m_{1})+(\alpha _{2}\circ f_{1})(l_{1}) & = & \alpha
_{2}(m_{1}^{\prime })+(\alpha _{2}\circ f_{1})(l_{1}^{\prime }) &  \\ 
\alpha _{2}(m_{1}^{\prime })+(f_{2}\circ \alpha _{1})(l_{1}) & = & \alpha
_{2}(m_{1}^{\prime })+(f_{2}\circ \alpha _{1})(l_{1}^{\prime }) &  \\ 
(f_{2}\circ \alpha _{1})(l_{1}) & = & (f_{2}\circ \alpha _{1})(l_{1}^{\prime
}) & \text{(}\alpha _{2}\text{ is cancellative)} \\ 
l_{1} & = & l_{1}^{\prime } & \text{(}f_{2}\text{ and }\alpha _{1}\text{ are
injective)} \\ 
m_{1}+f_{1}(l_{1}^{\prime }) & = & m_{1}^{\prime }+f_{1}(l_{1}^{\prime }) & 
\\ 
m_{1} & = & m_{1}^{\prime } & \text{(}f_{1}\text{ is cancellative)}%
\end{array}%
\end{equation*}
\end{enumerate}

\item We claim that $\alpha _{2}$ is a semi-epimorphism.

Let $m_{2}\in M_{2}.$ Since $\alpha _{3}$ and $g_{1}$ are surjective, there
exists $m_{1}\in M_{1}$ such that $g_{2}(m_{2})=(\alpha _{3}\circ
g_{1})(m_{1})=(g_{2}\circ \alpha _{2})(m_{1}).$ Since $g_{2}$ is $k$%
-uniform, $f_{2}(L_{2})=\mathrm{Ker}(g_{2})$ and $\alpha _{1}$ is
surjective, there exist $l_{1},l_{1}^{\prime }\in L_{1}$ such that%
\begin{eqnarray*}
m_{2}+(f_{2}\circ \alpha _{1})(l_{1}) &=&\alpha _{2}(m_{1})+(f_{2}\circ
\alpha _{1})(l_{1}^{\prime }) \\
m_{2}+\alpha _{2}(f_{1}(l_{1})) &=&\alpha _{2}(m_{1}+f_{1}(l_{1}^{\prime })).
\end{eqnarray*}%
Consequently, $M_{2}=\overline{\alpha _{2}(M_{1})},$ \emph{i.e.} $\alpha
_{2} $ is a semi-epimorphism. If $\alpha _{2}$ is $i$-uniform, then $M_{2}=%
\overline{\alpha _{2}(M_{1})}=\alpha _{2}(M_{1}),$ \emph{i.e.} $\alpha _{2}$
is surjective.$\blacksquare $
\end{enumerate}
\end{Beweis}

\begin{corollary}
\label{cor-short5}Consider the following commutative diagram of semimodules
with exact rows and assume that $M_{1}$ and $M_{2}$ are cancellative%
\begin{equation*}
\xymatrix{& L_1 \ar[r]^{f_1} \ar[d]_{\alpha_1} & M_1 \ar[r]^{g_1}
\ar[d]_{\alpha_2} & N_1 \ar[d]_{\alpha_3} \ar[r] & 0\\ 0 \ar[r] & L_2
\ar[r]^{f_2} & M_2 \ar[r]^{g_2} & N_2 & }
\end{equation*}

\begin{enumerate}
\item Let $\alpha _{2}$ be an isomorphism. Then $\alpha _{1}$ is surjective
if and only if $\alpha _{3}$ is injective.

\item Let $\alpha _{2}$ be $i$-uniform. If $\alpha _{1}$ and $\alpha _{3}$
are isomorphisms, then $\alpha _{2}$ is an isomorphism.
\end{enumerate}
\end{corollary}

\begin{proposition}
\label{short-5}\emph{(The Short Five Lemma)} Consider the following
commutative diagram of semimodules with $M_{1},M_{2}$ cancellative%
\begin{equation*}
\xymatrix{0 \ar[r] & L_1 \ar[r]^{f_1} \ar[d]_{\alpha_1} & M_1 \ar[r]^{g_1}
\ar[d]_{\alpha_2} & N_1 \ar[d]_{\alpha_3} \ar[r] & 0\\ 0 \ar[r] & L_2
\ar[r]^{f_2} & M_2 \ar[r]^{g_2} & N_2 \ar[r] & 0}
\end{equation*}%
Then $\alpha _{1},$ $\alpha _{3}$ are isomorphisms and $\alpha _{2}$ is $i$%
-uniform if and only if $\alpha _{2}$ is an isomorphism.
\end{proposition}

\begin{lemma}
\label{5-details}Consider the following commutative diagram of semimodules
with exact rows%
\begin{equation*}
\xymatrix{U_1 \ar[r]^{d_1} \ar[d]_{\gamma} & L_1 \ar[r]^{f_1}
\ar[d]_{\alpha_1} & M_1 \ar[r]^{g_1} \ar[d]_{\alpha_2} & N_1 \ar[r]
\ar[d]_{\alpha_3} \ar[r]^{h_1} & V_1 \ar[d]_{\delta} \\ U_2 \ar[r]^{d_2} &
L_2 \ar[r]^{f_2} & M_2 \ar[r]^{g_2} & N_2 \ar[r]^{h_2} & V_2}
\end{equation*}

\begin{enumerate}
\item Let $\gamma $ be surjective.

\begin{enumerate}
\item If $\alpha _{1}$ is injective and $\alpha _{3}$ is a
semi-monomorphisms, then $\alpha _{2}$ is a semi-monomorphism.

\item Assume that $f_{1}$ and $\alpha _{2}$ are cancellative. If $\alpha
_{1} $ and $\alpha _{3}$ are injective, then $\alpha _{2}$ is injective.
\end{enumerate}

\item Let $\delta $ be a semi-monomorphism. If $\alpha _{1},$ $\alpha _{3}$
are surjective (and $\alpha _{2}$ is $i$-uniform), then $\alpha _{2}$ is a
semi-epimorphism (surjective).

\item Let $f_{1},\alpha _{2}$ be cancellative, $\gamma $ be surjective and $%
\delta $ be injective. If $\alpha _{1}$ and $\alpha _{3}$ are isomorphisms,
then $\alpha _{2}$ is injective and a semi-epimorphism.
\end{enumerate}
\end{lemma}

\begin{Beweis}
Assume that the diagram is commutative and that the two rows are exact.

\begin{enumerate}
\item Let $\gamma $ be surjective.

\begin{enumerate}
\item We claim that $\alpha _{2}$ is a semi-monomorphism.

Suppose that $\alpha _{2}(m_{1})=0$ for some $m_{1}\in M_{1}.$ Then $(\alpha
_{3}\circ g_{1})(m_{1})=(g_{2}\circ \alpha _{2})(m_{1})=0.$ Since $\alpha
_{3}$ is a semi-monomorphism, $g_{1}(m_{1})=0$ and so $m_{1}=f_{1}(l_{1})$
for some $l_{1}\in L_{1}.$ It follows that $0=\alpha _{2}(m_{1})=(\alpha
_{2}\circ f_{1})(l_{1})=(f_{2}\circ \alpha _{1})(l_{1}),$ whence $\alpha
_{1}(l_{1})=(d_{2}\circ \gamma )(u_{1})=(\alpha _{1}\circ d_{1})(u_{1})$ for
some $u_{1}\in U_{1}$ (since $\gamma $ is surjective and $\mathrm{Ker}%
(f_{2})=d_{2}(U_{2})$). Since $\alpha _{1}$ is injective, $%
l_{1}=d_{1}(u_{1}) $ whence $m_{1}=f_{1}(l_{1})=(f_{1}\circ d_{1})(u_{1})=0.$

\item We claim that $\alpha _{2}$ is injective.

Suppose that $\alpha _{2}(m_{1})=\alpha _{2}(m_{1}^{\prime })$ for some $%
m_{1},m_{1}^{\prime }\in M_{1}.$ Then $(\alpha _{3}\circ
g_{1})(m_{1})=(g_{2}\circ \alpha _{2})(m_{1})=(g_{2}\circ \alpha
_{2})(m_{1}^{\prime })=(\alpha _{3}\circ g_{1})(m_{1}^{\prime }),$ whence $%
g_{1}(m_{1})=g_{1}(m_{1}^{\prime })$ because $\alpha _{3}$ is injective.
Since $g_{1}$ is $k$-uniform and $\mathrm{Ker}(g_{1})=f_{1}(L_{1}),$ there
exist $l_{1},l_{1}^{\prime }\in L_{1}$ such that $m_{1}+f_{1}(l_{1})=m_{1}^{%
\prime }+f_{1}(l_{1}^{\prime }).$ Then we have%
\begin{equation*}
\begin{array}{rclc}
\alpha _{2}(m_{1})+(\alpha _{2}\circ f_{1})(l_{1}) & = & \alpha
_{2}(m_{1}^{\prime })+(\alpha _{2}\circ f_{1})(l_{1}^{\prime }) &  \\ 
\alpha _{2}(m_{1}^{\prime })+(f_{2}\circ \alpha _{1})(l_{1}) & = & \alpha
_{2}(m_{1}^{\prime })+(f_{2}\circ \alpha _{1})(l_{1}^{\prime }) &  \\ 
f_{2}(\alpha _{1}(l_{1})) & = & f_{2}(\alpha _{1}(l_{1}^{\prime })) & \text{(%
}\alpha _{2}\text{ is cancellative)} \\ 
\alpha _{1}(l_{1})+k_{2} & = & \alpha _{1}(l_{1}^{\prime })+k_{2}^{\prime }
& \text{(}f_{2}\text{ is $k$-uniform)} \\ 
\alpha _{1}(l_{1})+(d_{2}\circ \gamma )(u_{1}) & = & \alpha
_{1}(l_{1}^{\prime })+(d_{2}\circ \gamma )(u_{1}^{\prime }) & \text{(}\gamma 
\text{ surjective, }\mathrm{Ker}(f_{2})=d_{2}(U_{2})\text{)} \\ 
\alpha _{1}(l_{1})+(\alpha _{1}\circ d_{1})(u_{1}) & = & \alpha
_{1}(l_{1}^{\prime })+(\alpha _{1}\circ d_{1})(u_{1}^{\prime }) &  \\ 
l_{1}+d_{1}(u_{1}) & = & l_{1}^{\prime }+d_{1}(u_{1}^{\prime }) & \text{(}%
\alpha _{1}\text{ is injective)} \\ 
f_{1}(l_{1}) & = & f_{1}(l_{1}^{\prime }) & \text{(}f_{1}\circ d_{1}=0\text{)%
} \\ 
m_{1}+f_{1}(l_{1}) & = & m_{1}+f_{1}(l_{1}^{\prime }) &  \\ 
m_{1}^{\prime }+f_{1}(l_{1}^{\prime }) & = & m_{1}+f_{1}(l_{1}^{\prime }) & 
\\ 
m_{1}^{\prime } & = & m_{1} & \text{(}f_{1}\text{ is cancellative)}%
\end{array}%
\end{equation*}
\end{enumerate}

\item Let $m_{2}\in M_{2}.$ Since $\alpha _{3}$ is surjective, there exists $%
n_{1}\in N_{1}$ such that $g_{2}(m_{2})=\alpha _{3}(n_{1}).$ It follows that 
$0=(h_{2}\circ g_{2})(m_{2})=(h_{2}\circ \alpha _{3})(n_{1})=(\delta \circ
h_{1})(n_{1}),$ whence $h_{1}(n_{1})=0$ since $\delta $ is a
semi-monomorphism. Since $g_{1}(M_{1})=\mathrm{Ker}(h_{1}),$ we have $%
n_{1}=g_{1}(m_{1})$ for some $m_{1}\in M_{1}.$ Notice that $(g_{2}\circ
\alpha _{2})(m_{1})=(\alpha _{3}\circ g_{1})(m_{1})=\alpha
_{3}(n_{1})=g_{2}(m_{2}).$ Since $g_{2}$ is $k$-uniform, $f_{2}(L_{2})=%
\mathrm{Ker}(g_{2})$ and $\alpha _{1}$ is surjective, there exists $%
l_{1},l_{1}^{\prime }\in L_{1}$ such that%
\begin{eqnarray*}
\alpha _{2}(m_{1})+(f_{2}\circ \alpha _{1})(l_{1}) &=&m_{2}+(f_{2}\circ
\alpha _{1})(l_{1}^{\prime }) \\
\alpha _{2}(m_{1}+f_{1}(l_{1})) &=&m_{2}+\alpha _{2}(f_{1}(l_{1}^{\prime })),
\end{eqnarray*}%
\emph{i.e.} $m_{2}\in \overline{\alpha _{2}(M_{1})}.$ Consequently, $M_{2}=%
\overline{\alpha _{2}(M_{1})}.$ If $\alpha _{2}$ is $i$-uniform, then $%
\alpha _{2}(M_{1})=\overline{\alpha _{2}(M_{1})}=M_{2},$ \emph{i.e. }$\alpha
_{2}$ is surjective.

\item This is a combination of \textquotedblleft 1\textquotedblright\ and
\textquotedblleft 2\textquotedblright .$\blacksquare $
\end{enumerate}
\end{Beweis}

\begin{corollary}
\label{5-lemma}\emph{(The Five Lemma) }Consider the following commutative
diagram of semimodules with exact rows and columns and assume that $M_{1}$
and $M_{2}$ are cancellative%
\begin{equation*}
\xymatrix{ & & & & 0 \ar[d] \\ U_1 \ar[r]^{e_1} \ar[d]_{\gamma} & L_1
\ar[r]^{f_1} \ar[d]_{\alpha_1} & M_1 \ar[r]^{g_1} \ar[d]_{\alpha_2} & N_1
\ar[r] \ar[d]_{\alpha_3} \ar[r]^{h_1} & V_1 \ar[d]_{\delta} \\ U_2
\ar[r]^{e_2} \ar[d] & L_2 \ar[r]^{f_2} & M_2 \ar[r]^{g_2} & N_2 \ar[r]^{h_2}
& V_2 \\ 0 & & & & }
\end{equation*}

\begin{enumerate}
\item If $\alpha _{1}$ and $\alpha _{3}$ are injective, then $\alpha _{2}$
is injective.

\item Let $\alpha _{2}$ be $i$-uniform. If $\alpha _{1}$ and $\alpha _{3}$
are surjective, then $\alpha _{2}$ is surjective.

\item Let $\alpha _{2}$ be $i$-uniform. If $\alpha _{1}$ and $\alpha _{3}$
are isomorphisms, then $\alpha _{2}$ is an isomorphism.
\end{enumerate}
\end{corollary}

\subsection*{The Nine Lemma}

\begin{lemma}
\label{9-1}Consider the following commutative diagram with exact columns and
assume that the second row is exact.%
\begin{equation*}
\xymatrix{ & & 0 \ar[d] & 0 \ar[d] & \\ & L_1 \ar[r]^{f_1} \ar[d]_{\alpha_1}
& M_1 \ar[r]^{g_1} \ar[d]_{\alpha_2} & N_1 \ar[d]_{\alpha_3} & \\ & L_2
\ar[r]^{f_2} \ar[d]_{\beta_1} & M_2 \ar[r]^{g_2} \ar[d]_{\beta_2} & N_2
\ar[d]_{\beta_3} &\\ & L_3 \ar[r]^{f_3} & M_3 \ar[r]^{g_3} & N_3 & }
\end{equation*}

\begin{enumerate}
\item If $f_{3}$ is injective and $f_{2}$ is cancellative, then the first
row is exact.

\item If $g_{2},$ $\beta _{1}$ are surjective, the third row is exact (and $%
g_{1}$ is $i$-uniform), then $g_{1}$ is a semi-epimorphism (surjective).
\end{enumerate}
\end{lemma}

\begin{Beweis}
Assume that the second row is exact.

\begin{enumerate}
\item Notice that $\alpha _{3}\circ g_{1}\circ f_{1}=g_{2}\circ \alpha
_{2}\circ f_{1}=g_{2}\circ f_{2}\circ \alpha _{1}=0,$ whence $g_{1}\circ
f_{1}=0$ since $\alpha _{3}$ is a monomorphism. In particular, $%
f_{1}(L_{1})\subseteq \mathrm{Ker}(g_{1}).$

\begin{itemize}
\item We claim that $\mathrm{Ker}(g_{1})\subseteq f_{1}(L_{1}).$

Let $m_{1}\in \mathrm{Ker}(g_{1}),$ \emph{i.e.} $g_{1}(m_{1})=0.$ It follows
that%
\begin{equation*}
\begin{array}{rclc}
(\alpha _{3}\circ g_{1})(m_{1}) & = & 0 &  \\ 
(g_{2}\circ \alpha _{2})(m_{1}) & = & 0 &  \\ 
\alpha _{2}(m_{1}) & = & f_{2}(l_{2}) & \text{(2nd row is proper exact)} \\ 
0 & = & (\beta _{2}\circ f_{2})(l_{2}) & \text{(}\beta _{2}\circ \alpha
_{2}=0\text{)} \\ 
0 & = & (f_{3}\circ \beta _{1})(l_{2}) &  \\ 
\beta _{1}(l_{2}) & = & 0 & \text{(}f_{3}\text{ is a semi-monomorphism)} \\ 
l_{2} & = & \alpha _{1}(l_{1}) & \text{(1st column is proper exact)} \\ 
f_{2}(l_{2}) & = & (f_{2}\circ \alpha _{1})(l_{1}) &  \\ 
\alpha _{2}(m_{1}) & = & \alpha _{2}(f_{1}(l_{1})) &  \\ 
m_{1} & = & f_{1}(l_{1}) & \text{(}\alpha _{2}\text{ is injective)}%
\end{array}%
\end{equation*}

\item We claim that $g_{1}$ is $k$-uniform.

Suppose that $g_{1}(m_{1})=g_{1}(m_{1}^{\prime })$ for some $%
m_{1},m_{1}^{\prime }\in M_{1}.$ It follows that%
\begin{equation*}
\begin{array}{rclc}
(\alpha _{3}\circ g_{1})(m_{1}) & = & (\alpha _{3}\circ g_{1})(m_{1}^{\prime
}) &  \\ 
(g_{2}\circ \alpha _{2})(m_{1}) & = & (g_{2}\circ \alpha _{2})(m_{1}^{\prime
}) &  \\ 
\alpha _{2}(m_{1})+f_{2}(l_{2}) & = & \alpha _{2}(m_{1}^{\prime
})+f_{2}(l_{2}^{\prime })\text{ (2nd row is exact)} &  \\ 
(\beta _{2}\circ f_{2})(l_{2}) & = & (\beta _{2}\circ f_{2})(l_{2}^{\prime })%
\text{ (}\beta _{2}\circ \alpha _{2}=0\text{)} &  \\ 
(f_{3}\circ \beta _{1})(l_{2}) & = & (f_{3}\circ \beta _{1})(l_{2}^{\prime })
&  \\ 
\beta _{1}(l_{2}) & = & \beta _{1}(l_{2}^{\prime })\text{ (}f_{3}\text{ is
injective)} &  \\ 
l_{2}+\alpha _{1}(l_{1}) & = & l_{2}^{\prime }+\alpha _{1}(l_{1}^{\prime })%
\text{ (first column is exact)} &  \\ 
f_{2}(l_{2})+(f_{2}\circ \alpha _{1})(l_{1}) & = & f_{2}(l_{2}^{\prime
})+(f_{2}\circ \alpha _{1})(l_{1}^{\prime }) &  \\ 
f_{2}(l_{2})+(\alpha _{2}\circ f_{1})(l_{1}) & = & f_{2}(l_{2}^{\prime
})+(\alpha _{2}\circ f_{1})(l_{1}^{\prime }) &  \\ 
\alpha _{2}(m_{1})+f_{2}(l_{2})+(\alpha _{2}\circ f_{1})(l_{1}) & = & \alpha
_{2}(m_{1})+f_{2}(l_{2}^{\prime })+(\alpha _{2}\circ f_{1})(l_{1}^{\prime })
&  \\ 
f_{2}(l_{2}^{\prime })+\alpha _{2}(m_{1}^{\prime }+f_{1}(l_{1})) & = & 
f_{2}(l_{2}^{\prime })+\alpha _{2}(m_{1}+f_{1}(l_{1}^{\prime })) &  \\ 
\alpha _{2}(m_{1}^{\prime }+f_{1}(l_{1})) & = & \alpha
_{2}(m_{1}+f_{1}(l_{1}^{\prime }))\text{ (}f_{2}\text{ is cancellative)} & 
\\ 
m_{1}^{\prime }+f_{1}(l_{1}) & = & m_{1}+f_{1}(l_{1}^{\prime })\text{ (}%
\alpha _{2}\text{ is injective)} & 
\end{array}%
\end{equation*}%
Since $f_{1}(L_{1})\subseteq \mathrm{Ker}(g_{1}),$ it follows that $g_{1}$
is $k$-uniform.
\end{itemize}

\item We claim that $g_{1}$ is a semi-epimorphism.

Let $n_{1}\in N_{1}$ and pick $m_{2}\in M_{2}$ such that $%
g_{2}(m_{2})=\alpha _{3}(n_{1})$ (by assumption $g_{2}$ is surjective). Then%
\begin{equation*}
\begin{array}{rclc}
g_{3}(\beta _{2}(m_{2})) & = & \beta _{3}(g_{2}(m_{2})) &  \\ 
& = & (\beta _{3}\circ \alpha _{3})(m_{2}) &  \\ 
& = & 0 & \text{(}\beta _{3}\circ \alpha _{3}=0\text{)} \\ 
\beta _{2}(m_{2}) & = & f_{3}(l_{3}) & \text{(3rd row is proper exact)} \\ 
& = & f_{3}(\beta _{1}(l_{2})) & \text{(}\beta _{1}\text{ is surjective)} \\ 
& = & \beta _{2}(f_{2}(l_{2})) &  \\ 
m_{2}+\alpha _{2}(m_{1}) & = & f_{2}(l_{2})+\alpha _{2}(m_{1}^{\prime }) & 
\text{(2nd column is exact)} \\ 
g_{2}(m_{2})+(g_{2}\circ \alpha _{2})(m_{1}) & = & (g_{2}\circ \alpha
_{2})(m_{1}^{\prime }) & \text{(}g_{2}\circ f_{2}=0\text{)} \\ 
\alpha _{3}(n_{1}+g_{1}(m_{1})) & = & \alpha _{3}(g_{1}(m_{1}^{\prime })) & 
\\ 
n_{1}+g_{1}(m_{1}) & = & g_{1}(m_{1}^{\prime }) & \text{(}\alpha _{3}\text{
is injective)}%
\end{array}%
\end{equation*}%
Consequently, $N_{1}=\overline{g_{1}(M_{1})}$ ($=$ $g_{1}(M_{1})$ if $g_{1}$
is assumed to be $i$-uniform).$\blacksquare $
\end{enumerate}
\end{Beweis}

\qquad Similarly, one can prove the following result:

\begin{lemma}
\label{9-3}Consider the following commutative diagram with exact columns and
assume that the second row is exact%
\begin{equation*}
\xymatrix{ & L_1 \ar[r]^{f_1} \ar[d]_{\alpha_1} & M_1 \ar[r]^{g_1}
\ar[d]_{\alpha_2} & N_1 \ar[d]_{\alpha_3} & \\& L_2 \ar[r]^{f_2}
\ar[d]_{\beta_1} & M_2 \ar[r]^{g_2} \ar[d]_{\beta_2} & N_2 \ar[d]_{\beta_3}
&\\ & L_3 \ar[r]^{f_3} \ar[d] & M_3 \ar[r]^{g_3} \ar[d] & N_3 & \\ & 0 & 0 &
}
\end{equation*}

\begin{enumerate}
\item If $g_{1}$ is surjective and $f_{3}$ is $i$-uniform, then the third
row is exact.

\item If $f_{2},$ $\alpha _{3}$ are injective, $\alpha _{2}$ is cancellative
and the first row is exact, then $f_{3}$ is injective.
\end{enumerate}
\end{lemma}

\begin{proposition}
\label{9}\emph{(The Nine Lemma) }Consider the following commutative diagram
with exact columns and assume that the second row is exact, $M_{2}$ is
cancellative and $f_{3},g_{1}$ are $i$-uniform%
\begin{equation*}
\xymatrix{ & 0 \ar@{.>}[d] & 0 \ar[d] & 0 \ar[d] & \\ 0 \ar@{.>}[r] & L_1
\ar[r]^{f_1} \ar[d]_{\alpha_1} & M_1 \ar[r]^{g_1} \ar[d]_{\alpha_2} & N_1
\ar[r] \ar[d]_{\alpha_3} & 0 \\ 0 \ar[r] & L_2 \ar[r]^{f_2} \ar[d]_{\beta_1}
& M_2 \ar[r]^{g_2} \ar[d]_{\beta_2} & N_2 \ar[r] \ar[d]_{\beta_3} & 0 \\ 0
\ar[r] & L_3 \ar[r]^{f_3} \ar[d] & M_3 \ar[r]^{g_3} \ar[d] & N_3
\ar@{-->}[r] \ar@{-->}[d] & 0 \\ & 0 & 0 & 0}
\end{equation*}%
Then the first row is exact if and only if the third row is exact.
\end{proposition}

\begin{Beweis}
The result follows immediately by combining Lemmas \ref{9-1} and \ref{9-3}.$%
\blacksquare $
\end{Beweis}

\subsection*{The Snake Lemma}

One of the basic homological lemmas that are proved usually in categories of
modules (\emph{e.g.} \cite{Wis1991}), or more generally in Abelian
categories, is the so called \emph{Kernel-Cokernel Lemma} (\emph{Snake Lemma}%
). Several versions of this lemma were proved also in non-abelian categories
(\emph{e.g.} \emph{homological categories} \cite{BB2004}, \emph{relative
homological categories} \cite{Jan2006} and incomplete relative homological
categories \cite{Jan2010b}).

\begin{theorem}
\label{snake}\emph{(The Snake Lemma) }Consider the following diagram of
semimodules in which the two middle squares are commutative and the two
middle rows are exact. Assume also that the columns are exact (or more
generally that $\alpha _{1},\alpha _{3}$ are $k$-uniform and $\alpha _{2}$
is uniform)%
\begin{equation*}
\xymatrix{ & 0 \ar[d] & 0 \ar[d] & 0 \ar[d] & \\ & {\rm Ker}(\alpha_1)
\ar[d]_{{\rm ker}(\alpha_1)} \ar@{.>}[r]^{f_K} & {\rm Ker}(\alpha_2)
\ar[d]_{{\rm ker}(\alpha_2)} \ar@{.>}[r]^{g_K} & {\rm Ker}(\alpha_3)
\ar[d]_{{\rm ker}(\alpha_3)} \ar@{-->}[dddll]^{\delta} & \\ & L_1
\ar[r]^{f_1} \ar[d]_{\alpha_1} & M_1 \ar[r]^{g_1} \ar[d]_{\alpha_2} & N_1
\ar[r] \ar[d]_{\alpha_3} & 0 \\ 0 \ar[r] & L_2 \ar[r]^{f_2} \ar[d]_{{\rm
coker}(\alpha_1)} & M_2 \ar[r]^{g_2} \ar[d]_{{\rm coker}(\alpha_2)} & N_2
\ar[d]_{{\rm coker}(\alpha_3)} & \\ & {\rm Coker}(\alpha_1)
\ar@{.>}[r]_{f_C} \ar[d] & {\rm Coker}(\alpha_2) \ar@{.>}[r]_{g_C} \ar[d] &
{\rm Coker}(\alpha_3) \ar[d] & \\ & 0 & 0 & 0}
\end{equation*}

\begin{enumerate}
\item There exist unique morphisms $f_{K},g_{K},f_{C}$ and $g_{C}$ which
extend the diagram commutatively.

\item If $f_{1}$ is cancellative, then the first row is exact.

\item If $f_{C}$ is $i$-uniform, then the last row is exact.

\item There exists a $k$-uniform \emph{connecting morphism} $\delta :\mathrm{%
Ker}(\alpha _{3})\longrightarrow \mathrm{Coker}(\alpha _{1})$ such that $%
\mathrm{Ker}(\delta )=\overline{g_{K}(\mathrm{Ker}(\alpha _{2}))}$ and $%
\delta (\mathrm{Ker}(\alpha _{3}))=\mathrm{Ker}(f_{C}).$

\item If $\alpha _{2}$ is cancellative and $g_{K}$ is $i$-uniform, then the
following sequence is exact%
\begin{equation*}
\xymatrix{ & {\rm Ker}(\alpha_2) \ar@{.>}[r]^{g_K} & {\rm Ker}(\alpha_3)
\ar@{-->}[r]^{\delta} & {\rm Coker}(\alpha_1) \ar@{.>}[r]^{f_C} & {\rm
Coker}(\alpha_2) &}
\end{equation*}
\end{enumerate}
\end{theorem}

\begin{Beweis}
\begin{enumerate}
\item The existence and uniqueness of the morphisms $f_{K},g_{K},f_{C}$ and $%
g_{C}$ is guaranteed by the definition of the (co)kernels and the
commutativity of the middle two squares.

\item This follows from Lemma \ref{9-1} applied to the first three rows.

\item This follows from Lemma \ref{9-3} applied to the last three rows.

\item We show first that $\delta $ exists and is well-defined.

\begin{itemize}
\item We define $\delta $ as follows: For $k_{3}\in \mathrm{Ker}(\alpha
_{3}),$ we choose $m_{1}\in M_{1}$ and $l_{2}\in L_{2}$ such that $%
g_{1}(m_{1})=k_{3}$ and $f_{2}(l_{2})=\alpha _{2}(m_{1});$ notice that this
is possible since $g_{1}$ is surjective and $(g_{2}\circ \alpha
_{2})(m_{1})=(\alpha _{3}\circ g_{1})(m_{1})=\alpha _{3}(k_{3})=0$ whence $%
\alpha _{2}(m_{1})\in \mathrm{Ker}(g_{2})=f_{2}(L_{2}).$ Define $\delta
(k_{3}):=\mathrm{coker}(\alpha _{1})(l_{2})=[l_{2}],$ the equivalence class
of $L_{2}/\alpha _{1}(L_{1})$ which contains $l_{2}.$

\item $\delta $ is well-defined, \emph{i.e.} $\delta (k_{3})$ is independent
of our choice of $m_{1}\in M_{1}$ and $l_{2}\in L_{2}$ satisfying the stated
conditions.

Suppose that $g_{1}(m_{1})=k_{3}=g_{1}(m_{1}^{\prime })$ for some $%
m_{1},m_{1}^{\prime }\in M_{1}$ and that $f_{2}(l_{2})=\alpha _{2}(m_{1}),$ $%
f_{2}(l_{2}^{\prime })=\alpha _{2}(m_{1}^{\prime })$ for some $%
l_{2},l_{2}^{\prime }\in L_{2}.$ Since the second row is exact, there exist $%
l_{1},l_{1}^{\prime }\in L_{1}$ such that $m_{1}+f_{1}(l_{1})=m_{1}^{\prime
}+f_{1}(l_{1}^{\prime }).$ It follows that%
\begin{equation*}
\begin{array}{rclc}
\alpha _{2}(m_{1})+(\alpha _{2}\circ f_{1})(l_{1}) & = & \alpha
_{2}(m_{1}^{\prime })+(\alpha _{2}\circ f_{1})(l_{1}^{\prime }) &  \\ 
f_{2}(l_{2})+(f_{2}\circ \alpha _{1})(l_{1}) & = & f_{2}(l_{2}^{\prime
})+(f_{2}\circ \alpha _{1})(l_{1}^{\prime }) &  \\ 
f_{2}(l_{2}+\alpha _{1}(l_{1})) & = & f_{2}(l_{2}^{\prime }+\alpha
_{1}(l_{1}^{\prime })) &  \\ 
l_{2}+\alpha _{1}(l_{1}) & = & l_{2}^{\prime }+\alpha _{1}(l_{1}^{\prime })
& \text{(}f_{2}\text{ is injective)} \\ 
\lbrack l_{2}] & = & [l_{2}^{\prime }] & 
\end{array}%
\end{equation*}%
Thus $l_{2}$ and $l_{2}^{\prime }$ lie in the same equivalence class of $%
L_{2}/\alpha _{1}(L_{1}),$ \emph{i.e.} $\delta $ is well-defined.

\item Clearly $\overline{g_{K}(\mathrm{Ker}(\alpha _{2}))}\subseteq \mathrm{%
Ker}(\delta )$ (notice that $f_{2}$ is a semi-monomorphism). We claim that $%
\mathrm{Ker}(\delta )\subseteq \overline{g_{K}(\mathrm{Ker}(\alpha _{2}))}.$

Let $k_{3}\in \mathrm{Ker}(\delta )$ and pick some $m_{1}\in M_{1}$ and $%
l_{2}\in L_{2}$ such that $g_{1}(m_{1})=k_{3}$ and $f_{2}(l_{2})=\alpha
_{2}(m_{1}).$ By assumption, $[l_{2}]=\delta (k_{3})=0,$ \emph{i.e.} $%
l_{2}+\alpha _{1}(l_{1})=\alpha _{1}(l_{1}^{\prime })$ for some $%
l_{1},l_{1}^{\prime }\in L_{1}.$Then we have%
\begin{equation*}
\begin{array}{rclc}
f_{2}(l_{2})+(f_{2}\circ \alpha _{1})(l_{1}) & = & (f_{2}\circ \alpha
_{1})(l_{1}^{\prime }) &  \\ 
\alpha _{2}(m_{1})+\alpha _{2}(f_{1}(l_{1})) & = & \alpha
_{2}(f_{1}(l_{1}^{\prime })) &  \\ 
m_{1}+f_{1}(l_{1})+k_{2} & = & f_{1}(l_{1}^{\prime })+k_{2}^{\prime } & 
\text{(}\alpha _{2}\text{ is }k\text{-uniform)} \\ 
k_{3}+g_{K}(k_{2}) & = & g_{K}(k_{2}^{\prime }) & \text{(}g_{1}\circ f_{1}=0%
\text{)}%
\end{array}%
\end{equation*}%
Consequently, $\overline{g_{K}(\mathrm{Ker}(\alpha _{2}))}=\mathrm{Ker}%
(\delta ).$

\item Let $k_{3}\in \mathrm{Ker}(\alpha _{3})$ and pick some $m_{1}\in
M_{1}, $ $l_{2}\in L_{2}$ such that $g_{1}(m_{1})=k_{3}$ and $%
f_{2}(l_{2})=\alpha _{2}(m_{1}).$ Then we have 
\begin{equation*}
(f_{C}\circ \delta )(k_{3})=f_{C}([l_{2}])=[f_{2}(l_{2})]=[\alpha
_{2}(m_{1})]=[0].
\end{equation*}%
Consequently, $\delta (\mathrm{Ker}(\alpha _{3}))\subseteq \mathrm{Ker}%
(f_{C}).$ We claim that $\mathrm{Ker}(f_{C})\subseteq \delta (\mathrm{Ker}%
(\alpha _{3})).$

Let $[l_{2}]\in \mathrm{Ker}(f_{C}),$ so that $%
[f_{2}(l_{2})]=f_{C}([l_{2}])=[0].$ Then there exist $m_{1},m_{1}^{\prime
}\in M_{1}$ such that $f_{2}(l_{2})+\alpha _{2}(m_{1})=\alpha
_{2}(m_{1}^{\prime }).$ Since $\alpha _{2}$ is $i$-uniform, there exists $%
\mathbf{m}_{1}\in M_{1}$ such that $\alpha _{2}(\mathbf{m}%
_{1})=f_{2}(l_{2}). $ It follows that $(\alpha _{3}\circ g_{1})(\mathbf{m}%
_{1})=(g_{2}\circ \alpha _{2})(\mathbf{m}_{1})=(g_{2}\circ f_{2})(l_{2})=0.$
So, $g_{1}(\mathbf{m}_{1})\in \mathrm{Ker}(\alpha _{3})$ and $[l_{2}]=\delta
(g_{1}(\mathbf{m}_{1})).$ Consequently, $\mathrm{Ker}(f_{C})=\delta (\mathrm{%
Ker}(\alpha _{3})).$

\item We claim that $\delta $ is $k$-uniform.

Suppose that $\delta (k_{3})=\delta (k_{3}^{\prime })$ for some $%
k_{3},k_{3}^{\prime }\in \mathrm{Ker}(\alpha _{3})$ and pick $%
m_{1},m_{1}^{\prime }\in M_{1},$ $l_{2},l_{2}^{\prime }\in L_{2}$ such that $%
g_{1}(m_{1})=k_{3},$ $g_{1}(m_{1}^{\prime })=k_{3}^{\prime },$ $\alpha
_{2}(m_{1})=f_{2}(l_{2})$ and $\alpha _{2}(m_{1}^{\prime
})=f_{2}(l_{2}^{\prime }).$ By assumption, $[l_{2}]=[l_{2}^{\prime }],$ 
\emph{i.e. }$l_{2}+\alpha _{1}(l_{1})=l_{2}^{\prime }+\alpha
_{1}(l_{1}^{\prime })$ for some $l_{1},l_{1}^{\prime }\in L_{1}.$ It follows
that%
\begin{equation*}
\begin{array}{rclc}
f_{2}(l_{2})+(f_{2}\circ \alpha _{1})(l_{1}) & = & f_{2}(l_{2}^{\prime
})+(f_{2}\circ \alpha _{1})(l_{1}^{\prime }) &  \\ 
\alpha _{2}(m_{1})+(\alpha _{2}\circ f_{1})(l_{1}) & = & \alpha
_{2}(m_{2}^{\prime })+(\alpha _{2}\circ f_{1})(l_{1}^{\prime }) &  \\ 
m_{1}+f_{1}(l_{1})+k_{2} & = & m_{1}^{\prime }+f_{1}(l_{1}^{\prime
})+k_{2}^{\prime } & \text{(}\alpha _{2}\text{ is }k\text{-uniform)} \\ 
g_{1}(m_{1})+g_{K}(k_{2}) & = & g_{1}(m_{1}^{\prime })+g_{K}(k_{2}^{\prime })
& \text{(}g_{1}\circ f_{1}=0\text{)} \\ 
k_{3}+g_{K}(k_{2}) & = & k_{3}^{\prime }+g_{K}(k_{2}^{\prime }) & 
\end{array}%
\end{equation*}%
\qquad Since $g_{K}(\mathrm{Ker}(\alpha _{2}))\subseteq \mathrm{Ker}(\delta
),$ we conclude that $\delta $ is $k$-uniform.
\end{itemize}

\item If $g_{K}$ is $i$-uniform, then we have $\mathrm{Ker}(\delta )=%
\overline{g_{K}(\mathrm{Ker}(\alpha _{2}))}=g_{K}(\mathrm{Ker}(\alpha _{2}))$
and it remains only to prove that $f_{C}$ is $k$-uniform.

Suppose that $f_{C}[l_{2}]=f_{C}[l_{2}^{\prime }]$ for some $%
l_{2},l_{2}^{\prime }\in L_{2}.$ Then there exist $m_{1},m_{1}^{\prime }\in
M_{1}$ such that $f_{2}(l_{2})+\alpha _{2}(m_{1})=f_{2}(l_{2}^{\prime
})+\alpha _{2}(m_{1}^{\prime }).$ It follows that%
\begin{equation*}
\begin{array}{rclc}
(g_{2}\circ \alpha _{2})(m_{1}) & = & (g_{2}\circ \alpha _{2})(m_{1}^{\prime
})\text{ (}g_{2}\circ f_{2}=0\text{)} &  \\ 
(\alpha _{3}\circ g_{1})(m_{1}) & = & (\alpha _{3}\circ g_{1})(m_{1}^{\prime
}) &  \\ 
g_{1}(m_{1})+k_{3} & = & g_{1}(m_{1}^{\prime })+k_{3}^{\prime }\text{ (}%
\alpha _{3}\text{ is }k\text{-uniform)} &  \\ 
g_{1}(m_{1}+\mathbf{m}_{1}) & = & g_{1}(m_{1}^{\prime }+\mathbf{m}%
_{1}^{\prime })\text{ (}g_{1}\text{ is surjective)} &  \\ 
m_{1}+\mathbf{m}_{1}+f_{1}(\mathbf{l}_{1}) & = & m_{1}^{\prime }+\mathbf{m}%
_{1}^{\prime }+f_{1}(\mathbf{l}_{1}^{\prime })\text{ (2nd row is exact)} & 
\\ 
\alpha _{2}(m_{1})+\alpha _{2}(\mathbf{m}_{1})+(\alpha _{2}\circ f_{1})(%
\mathbf{l}_{1}) & = & \alpha _{2}(m_{1}^{\prime })+\alpha _{2}(\mathbf{m}%
_{1}^{\prime })+(\alpha _{2}\circ f_{1})(\mathbf{l}_{1}^{\prime }) &  \\ 
f_{2}(l_{2}^{\prime })+\alpha _{2}(m_{1})+\alpha _{2}(\mathbf{m}%
_{1})+(f_{2}\circ \alpha _{1})(\mathbf{l}_{1}) & = & [f_{2}(l_{2}^{\prime
})+\alpha _{2}(m_{1}^{\prime })]+\alpha _{2}(\mathbf{m}_{1}^{\prime
})+(f_{2}\circ \alpha _{1})(\mathbf{l}_{1}^{\prime }) &  \\ 
f_{2}(l_{2}^{\prime })+\alpha _{2}(m_{1})+\alpha _{2}(\mathbf{m}%
_{1})+(f_{2}\circ \alpha _{1})(\mathbf{l}_{1}) & = & f_{2}(l_{2})+\alpha
_{2}(m_{1})+\alpha _{2}(\mathbf{m}_{1}^{\prime })+(f_{2}\circ \alpha _{1})(%
\mathbf{l}_{1}^{\prime }) &  \\ 
f_{2}(l_{2}^{\prime })+\alpha _{2}(\mathbf{m}_{1})+(f_{2}\circ \alpha _{1})(%
\mathbf{l}_{1}) & = & f_{2}(l_{2})+\alpha _{2}(\mathbf{m}_{1}^{\prime
})+(f_{2}\circ \alpha _{1})(\mathbf{l}_{1}^{\prime })\text{ (}\alpha _{2}%
\text{ is cancellative)} &  \\ 
f_{2}(l_{2}^{\prime }+\mathbf{l}_{2}+\alpha _{1}(\mathbf{l}_{1})) & = & 
f_{2}(l_{2}+\mathbf{l}_{2}^{\prime }+\alpha _{1}(\mathbf{l}_{1}^{\prime }))%
\text{ (}g_{1}(\mathbf{m}_{1}),\text{ }g_{1}(\mathbf{m}_{1}^{\prime })\in 
\mathrm{Ker}(\alpha _{3})\text{)} &  \\ 
l_{2}^{\prime }+\mathbf{l}_{2}+\alpha _{1}(\mathbf{l}_{1}) & = & l_{2}+%
\mathbf{l}_{2}^{\prime }+\alpha _{1}(\mathbf{l}_{1}^{\prime })\text{ (}f_{2}%
\text{ is injective)} &  \\ 
\lbrack l_{2}^{\prime }]+[\mathbf{l}_{2}] & = & [l_{2}]+[\mathbf{l}%
_{2}^{\prime }] &  \\ 
\lbrack l_{2}^{\prime }]+\delta (k_{3}) & = & [l_{2}]+\delta (k_{3}^{\prime
})\text{ (definition of }\delta \text{)} & 
\end{array}%
\end{equation*}%
Since $\delta (\mathrm{Ker}(\alpha _{3}))\subseteq \mathrm{Ker}(f_{C}),$ we
conclude that $f_{C}$ is $k$-uniform.$\blacksquare $
\end{enumerate}
\end{Beweis}

\textbf{Acknowledgments.} The author thanks the anonymous referee for
her/his careful reading of the paper and the useful suggestions which
improved the revised version of this manuscript. He also thanks all
mathematicians who clarified to him some issues related to the nature of the
categories of semimodules and exact sequences or sent to him related
manuscripts especially G. Janelidze, Y. Katsov, F. Linton, A. Patchkoria, H.
Porst and R. Wisbauer.

\end{document}